\input amstex
\input amsppt.sty
\magnification\magstep1

\def\ni\noindent
\def\sbs{\subset}

\def\as{\operatorname{asdim}}
\def\diam{\operatorname{diam}}

\def\R{\text{\bf R}}

\def\Q{\text{\bf Q}}
\def\Z{\text{\bf Z}}

\def\N{\text{\bf N}}
\def\L{\text{\bf L}}

\def\sC{\Cal C}

\def\sU{\Cal U}

\hoffset= 0.0in
\voffset= 0.0in
\hsize=32pc
\vsize=38pc
\baselineskip=24pt
\NoBlackBoxes
\topmatter
\author
A.N. Dranishnikov
\endauthor

\title
HyperEuclidean manifolds and the Novikov Conjecture
\endtitle
\abstract
We develop some basic Lipschitz homotopy technique and apply it to
manifolds with finite asymptotic dimension. In particular we show that
the Higson compactification of a uniformly contractible manifold is mod $p$
acyclic in the finite dimensional case. Then we give an alternative proof
of the Higher Signature Novikov Conjecture for the groups with finite
asymptotic dimension.
Finally we define an asymptotically piecewise Euclidean metric space
as a space which admits an approximation by Euclidean asymptotic polyhedra.
We show that the Gromov-Lawson conjecture holds for the asymptotically
piecewise
Euclidean groups. Also we prove that expanders are not asymptotically
piecewise Euclidean
\endabstract

\thanks The author was partially supported by NSF grant DMS-9971709.
\endthanks

\address University of Florida, Department of Mathematics, P.O.~Box~118105,
358 Little Hall, Gainesville, FL 32611-8105, USA
\endaddress

\subjclass Primary 20H15
\endsubjclass

\email  dranish\@math.ufl.edu
\endemail

\keywords  dimension, asymptotic dimension,
absolute extensor, Higson corona, Novikov conjecture.
\endkeywords
\endtopmatter

\document
\head \S1 Introduction \endhead

The Novikov Conjecture states that the higher signatures of a manifold are
homotopy invariants. The ordinary signature of a manifold $M$ can be computed by
the Hirzebruch formula $\sigma(M)=\langle L,[M]\rangle$ where the
right hand side is the integration of the Hirzebruch polynomial over the
fundamental homology class of $M$. If the fundamental group 
$\Gamma=\pi_1(M)\ne 0$,
then
there are distinguished homology classes $\beta\in H_*(M)$ which come from
the group $\Gamma$. The integration of the Hirzebruch polynomial over these classes
$\beta$ gives rise the higher signatures $\sigma_{\beta}(M)$ of $M$. In view of
this it makes sense to speak about the Novikov Conjecture for a given finitely
presented group $\Gamma$.

The most successful approach to the Novikov Conjecture in the last decade was
the so called  coarse approach. The coarse approach consists of studying the
large scale geometry of a group $\Gamma$ as a metric space with the word metric.
Quite often it is more convenient to consider the universal cover $\tilde M$
with the induced metric, though in the coarse sense $\tilde M$ is equivalent to
$\Gamma$. The most advanced result on the Novikov Conjecture is
a theorem of G. Yu [Yu2] which asserts that the Novikov Conjecture holds true
for groups $\Gamma$ that admit an embedding in a coarse sense in the Hilbert
space $l_2$. This Yu's theorem formally generalizes the other his theorem [Yu1]
saying that
the Novikov Conjecture holds for groups $\Gamma$ with finite asymptotic dimension
$\as\Gamma$ . Both Yu's result are obtained by proving the
coarse Baum-Connes conjecture which states that the Roe index map
$$
A:K_*(\tilde M)\to K_*(C^*(\tilde M))
$$
is an isomorphism when $\tilde M$ is the universal cover of a closed aspherical
manifold $M$ [H-R],[Ro2]. The $C^*$-algebra $C^*(\tilde M)$ here is the completion of
the algebra of infinite matrices $A=(a_{x,y})_{x,y\in\Gamma}$ whose entries are
compact operators on $l_2$ and $a_{x,y}=0$ for $d_{\Gamma}(x,y)\ge r$ for some $r=r(A)$.
Thus both Yu's works are heavily operator algebra theoretic. The descent from
the coarse Baum-Connes conjecture to the Novikov Conjecture is shown in [Ro2].
It uses the homotopy fixed point theory and it proves the analytic Novikov
conjecture. The analytic Novikov conjecture states that the analytic assembly map
$A_{\Gamma}:K_*(B\Gamma)\to K_*(C^*_r\Gamma)$ is injective. Here $C^*_r\Gamma$ is
the reduced group $C^*$-algebra. The classical Novikov Conjecture is equivalent to
the statement that the rational assembly map from the surgery exact sequence
$H_*(B\Gamma;\L)\otimes\Q\to L_*(\Z[\Gamma])\otimes\Q$ is injective [K-M], [F-R-R].
It is known that the analytic
Novikov conjecture implies the original conjecture [F-R-R]. 

Despite on recent counterexamples to the coarse Baum-Connes conjecture,
this approach to the Novikov Conjecture is not exhausted yet.
It is possible to show
that the equivariant
split injectivity of the Roe index map also enables to derive the Novikov Conjecture.
As far as it known to the author,
the monomorphism version of the coarse Baum-Connes conjecture is not disproved
to this moment.

A more geometric coarse reduction of the Novikov Conjecture was considered by
Ferry and Weinberger [F-W]. They were looking for obstructions to make a homotopy
between two closed aspherical manifold tangential. This approach can be traced
back to Farrell-Hsiang [F-H]. The Ferry-Weinberger approach reduces the Novikov
Conjecture to the question whether the boundary homomorphism
$\delta:\check H^{n-1}(\nu\tilde M;\Q)\to H^n_c(\tilde M;\Q)$ is an equivariant
split epimorphism where $\nu\tilde M$ is the Higson corona [Ro1].

We recall that the Higson compactification $X\cup\nu X$ of $X$ is generated by
the algebra of functions with the gradient tending to zero at infinity.
The exact sequence of pair $(\tilde M\cup\nu\tilde M,\nu\tilde M)$ implies that
the boundary homomorphism 
$\delta$ is an epimorphism provided $H^n(\tilde M\cup\nu\tilde M;\Q)=0$.
The assertion that the Higson compactification  $\tilde M\cup\nu\tilde M$ is acyclic
is called the Higson conjecture [Ro1]. It is known that for integral coefficients
this conjecture false even for $\R^n$ [Ke],[D-F]. In \S 3 we prove the Higson
conjecture for finite coefficients in the finite dimensional case. There
are chances that the mod $p$ Higson conjecture holds in full generality.
In \S 5 we formulate a stable version of the Higson conjecture
which cannot be disproved by simple examples.

The other, in some sense equivalent, approach to the Novikov Conjecture
is due to Gromov which goes back to his work on the positive scalar curvature
[G-L].
He introduced the notions of hyperspherical and hypereuclidean manifolds
[G2]. Roe proved [Ro1] that an $n$-manifold $\tilde M$ is hypereuclidean if and only
if
the Higson corona $\nu\tilde M$ admits a map of degree one onto an $n-1$-sphere
$S^{n-1}$. It is easy to show that if a group $\Gamma$ acts on $\tilde M$ by
isometries then that action can be extended to an action on the Higson corona
[Dr1]. We introduce the notion of an equivariantly hypereuclidean 
(or $\Gamma$-hypereuclidean) $n$-manifold
as a manifold $\tilde M$ with a group $\Gamma$ acting on it properly and 
cocompactly whose Higson corona admits a map $f$ of
 degree one onto an $(n-1)$-dimensional sphere
such that the action of $\Gamma$ is fixed on $f^*(H^{n-1}(S^{n-1})$. 
Using Ferry-Weinberger
approach it is possible to show that the Novikov Conjecture holds for
manifolds $M$ with the equivariantly
$\pi_1(M)$-hypereuclidean universal cover $\tilde M$. In $\S 4$, $\S 5$ 
we show that if
$\as\pi_1(M)<\infty$ then $\tilde M\times\R^n$ is equivariantly hypereuclidean
for some $n$. This allows to establish the Novikov Conjecture for $M$.

The concepts of
hypersphericity is formally weaker. There is an open question whether an integrally
hyperspherical manifold is always hypereuclidean. The hypersphericity of
$\tilde M$ implies the Gromov-Lawson conjecture for $M$ [G-L]:
{\it A closed aspherical manifold cannot carry a metric of a positive scalar
curvature}.
This conjecture is quite close to the Novikov Conjecture
[Ros]. Perhaps the
Novikov Conjecture can be derived from some equivariant version of the
hypersphericity. In [Dr2] it was shown that in finite dimensional case
the manifold $\tilde M\times\R^k$ is hyperspherical for some $k$.

 In this paper in \S 6
we introduce the notion of {\it asymptotically piecewise Euclidean metric space} in the
coarse category
as a space that admits an approximation by piecewise Euclidean asymptotic polyhedra.
Then we
extend Yu's result about Gromov-Lawson conjecture to asymptotically 
piecewise Euclidean
manifolds $\tilde M$. Note that all groups uniformly embeddable in the Hilbert
space
$l_2$ are asymptotically piecewise Euclidean, although it is unclear whether the
inclusion of these classes is strong.  
 Like in the case with Yu's theorem [G4], this result is
also limited by expanders. This is demonstrated in \S 7, where
 we show that an asymptotically piecewise
Euclidean metric space cannot contain an expander.

We note that there is a
 connection between Ferry-Weinberger's (and Gromov's) 
and the coarse Baum-Connes
conjecture approaches  to the
Novikov Conjecture which is based on the fact that the topological K-theory
$K_*(\nu\tilde M)$
is an approximation
to $C^*$-algebra K-theory $K_*(C^*(\tilde M))$ [Ro1].

Using Roe's coarse cohomology one can define an asymptotic cohomological
dimension $\as_{\Z}X$ in the coarse category as it was done in [Dr1].
It seems to do that properly one has to make a shift in the grading of Roe's cohomology
(it was not done in [Dr1]).  Only in that case we would get the equality
$\as_{\Z}\R^n=n$. The Roe cohomology is dual in the macro-micro sense to
the Steenrod homology.  Curiously, Steenrod defined his homology also with shifted
dimensions [St] and only later Sitnikov made the correction [Sit].
We note that
for the universal covers $\tilde M$ of aspherical manifolds 
 always there is the
inequality $\as_{\Z}\tilde M<\infty$ [Dr1]. If the asymptotic cohomological dimension
 agreed with Gromov's asymptotic dimension, then
the Novikov Conjecture would follow. In the ordinary topology the problem
about coincidence of the cohomological dimension and the Lebesgue dimension was known
from late 20s as the Alexandroff problem. The Alexandroff problem was solved by
a counterexample [Dr3]. In the large scale world  a counterexample
was constructed in [D-F-W] but that example does not have a bounded geometry.
Recently a counterexample among finitely presented groups was found by Gromov
by means
of expanders and random groups [G].

The author wish to thank the Max Plank Institute fur Mathematik for the
hospitality where this paper was written.

\head \S2 Lipschitz homotopy \endhead

In this section we consider variations of the following question:
{\it When does a null
homotopic $\lambda$-Lipschitz map $f:X\to Y$ admit a $\mu$-Lipschitz homotopy
$H:X\times I\to Y$ to a constant map?}
We recall that a map $f:X\to Y$ between metric spaces is $\lambda$-Lipschitz
if $d_Y(f(x),f(x'))\le\lambda d_X(x,x')$ for all $x,x'\in X$.
Denote by $L(f)=\sup\{\frac{d_Y(f(x),f(x'))}{d_X(x,x')}\}$. Then for
a $\lambda$-Lipschitz map $f$ we have $L(f)\le\lambda$.

First we give an answer to the question for finite simplicial complexes.
Every simplicial complex $K$ carries a metric
where
all simplexes are the standard of size one.
We will call such metric {\it uniform} and usually we will
denote corresponding metric space as $K_U$. If it is not specified,
we will assume that a finite complex always supplied with the
uniform metric.

\proclaim{Lemma 2.1}
Suppose that $X$ and $Y$ are finite simplicial complexes.
Then for every $\lambda$
there exists $\mu=\mu(\lambda)$ such that every null homotopic
$\lambda$-Lipschitz map $f:X\to Y$ admits a $\mu$-Lipschitz homotopy
$H:X\times I\to Y$ to a constant map.
\endproclaim
\demo{Proof}
Let $r_n$ denote radius of the inscribed sphere in the standard $n$-simplex.
We fix a subdivision $T$ of $X$ with the mesh $<\frac{r_n}{4\lambda}$.
There are finitely many different simplicial maps $\phi:T\to Y$.
We consider only null homotopic maps $\phi_i$. For every $i$
we fix a homotopy $H_i:T\times I\to Y$ to a constant map. By 
Simplicial Approximation Theorem we may assume that the map $H_i$ is
$\mu_i$-Lipschitz. We take $\mu\ge\max\{2\mu_i\}$. Also we assume that
$\mu\ge\frac{1}{2r_T}$ where $r_T$ is the minimum of radii of inscribed
spheres in simplices of $T$. According to the following Lemma 2.2 every
$\lambda$-Lipschitz
map $f:X\to Y$ is homotopic to a simplicial map $g:T\to Y$ with respect to
$T$. Moreover, the corresponding homotopy is 
$(\mu/2)$-Lipschitz with the above $\mu$. Then $g=\phi_i$
for some $i$ and a $\mu$-Lipschitz homotopy of $f$ to a constant map will be
the sum of the above homotopy deforming $f$ to $g$ and $H_i$.\qed
\enddemo

By $B_r(x_0)=\{x\in X\mid d(x,x_0)\le r\}$ we denote a closed ball of
 radius $r$ in a metric space $X$ with center at $x_0\in X$.
The unit $n$-ball in $\R^n$ will be denoted as $B^n$.
\proclaim{Lemma 2.2}
Let $f:L\to K$ be a $\lambda$-Lipschitz map between uniform finite
dimensional polyhedra.
Then $f$ is homotopic to a simplicial map $g$ with respect to some
subdivision of $L$ by means of $\mu$-Lipschitz homotopy, where
$\mu$ depends on $\lambda$ and
$\dim L,\dim K$ only.
\endproclaim
\demo{Proof}
We apply the standard argument of the Simplicial Approximation theorem.
Consider an open cover $\sU=\{f^{-1}(OSt(v,K))\mid v\in K^{(0)}\}$ of $L$
where $OSt(v,K)$ denotes the open star of a vertex $v$ in a complex $K$.
Since $f$ is $\lambda$-Lipschitz, the Lebesgue number of $\sU$ is greater
than $\frac{r_n}{2\lambda}$, where $r_n$ denotes radius of the inscribed sphere in
the unit standard $n$-simplex and $n=\dim K$. To see that, take a point $y\in L$
and consider a closest vertex $v\in K^{(0)}$ to $f(y)$. Then the ball
$B_{\frac{r_n}{2}}(y)$ is contained in $OSt(v,K)$. Hence the ball
$B_{\frac{r_n}{2\lambda}}$ is contained in $f^{-1}(OSt(v,K))$. Consider a
triangulation $T$ of $L$ which is an iterated barycentric subdivision
of $K$ with mesh $<\frac{r_n}{4\lambda}$. Let $r_T$ denote the minimum of radii of
inscribed spheres in simplices of $T$. We take $\mu\ge\frac{1}{r_T}$. Then
$\mu$ depends on $\lambda, n$ and dimension of
$L$. We consider a simplicial approximation $g:T\to K$ defined by the
standard rule: for every vertex $a\in T^{(0)}$ we take $g(a)\in K^{(0)}$ such that
 $OSt(a,L)\subset f^{-1}(OSt(g(a),K))$.  Let us show that for every $x\in L$
the image $g(x)$ belongs to $\sigma$ whenever
$f(x)$ lies in the interior of a simplex $\sigma$ in $K$. It suffices to show
that $g(\Delta)=\sigma$, where
$\Delta$ is a unique simplex in $L$ containing $x$ as an interior point.
Let $p\in \Delta^{(0)}$ be a vertex, then $x\in OSt(p,L)$. Therefore
$f(x)\in f(OSt(p,L))\subset OSt(g(p),K)$. This means that $g(p)$ is a
vertex of $\sigma$. Then a Lipschitz homotopy between $f$ and $g$ is defined
linearly by joining $f(x)$ with $g(x)$ in $\sigma$. \qed
\enddemo

REMARK. Lemma 2.1 holds true if one consider a compact metric space $X$.
In that case by the Ascoli-Arzela theorem the space  $Y^X_{\lambda}$ of
$\lambda$-Lipschitz maps $f:X\to Y$ is compact. We can take a finite
$\epsilon$-net in  $Y^X_{\lambda}$ for small enough $\epsilon$.
Similarly, for every map $\phi_i:X\to Y$ from this $\epsilon$-net we fix
a Lipschitz homotopy $H_i$ to a constant map. Then any other
$\lambda$-Lipschitz map is $\mu'$-homotopic to one from the net.

The following lemma can be derived also from results of Siegel and Williams [S-W].
\proclaim{Lemma 2.3}
Let $Y$ be a finite simplicial complex with $\pi_n(Y)$ finite. Then for every $\lambda$
there is $\mu$ such that every map $f:B^n\to Y$ with $L(f|_{S^{n-1}})\le\lambda$
can be deformed to a $\mu$-Lipschitz map $g:B^n\to Y$ by means of a homotopy
$h_t:B^n\to Y$ with $h_t|_{S^{n-1}}=f|_{S^{n-1}}$.
\endproclaim
\demo{Proof}
We consider a finite family of simplicial maps $\phi_i:S^{n-1}\to Y$ as in
the proof of Lemma 2.1. Since the group $\pi_n(Y)$ is finite, there are finite
number of homotopically different extensions $H^j_i:B^n\to Y$.
We assume that every $H^j_i$ is $\mu^j_i$-Lipschitz. The rest of the proof
is the same as in Lemma 1.\qed
\enddemo
\proclaim{Lemma 2.4}
Let $L$ be a finite dimensional complex and
let $K$ be a finite complex with finite homotopy groups
$\pi_i(K)$ for $i\le\dim L+1$. Let
$f,g:L\to K$ be homotopic Lipschitz maps. Then every homotopy between $f$
and $g$ can be deformed to a Lipschitz homotopy
$H:L\times[0,1]\to K$.
\endproclaim
\demo{Proof}
Let $F:L\times I\to K$ be a homotopy between $f$ and $g$.
By induction on $n$,  and using Lemma 2.3 we construct a $\mu_n$-Lipschitz map
$H_n:L^{(n)}\times I\cup L\times\{0,1\}\to K$
which is a deformation of $F$ restricted to the $n$-skeleton $L^{(n)}$
such that $H_n$ extends $H_{n-1}$. Then $H=H_m$ for $m=\dim L$.
\qed
\enddemo

\

\head \S3 Modulo $p$ Higson conjecture\endhead

A countable simplicial complex $L$ with a geodesic metric on it is called
an {\it asymptotic polyhedron} if
every simplex in $L$ is isometric to a simplex $\Delta$ spanned in the
Hilbert space and
the radii of inscribed spheres $r_{\Delta}$ tend to infinity.
In this definition any
Banach space can be used. Moreover,
asymptotic polyhedra naturally appear with $l_{\infty}$ metric.
Since we are working in this section with
finite dimensional complexes we can consider only the Euclidean case.

We recall that by the definition the {\it asymptotic dimension} $\as X$
of a metric space $X$ does not exceed $n$ if for any arbitrary large number $d$
there is a uniformly bounded open cover $\sU$ of $X$ with multiplicity $\le n+1$
and with the Lebesgue number $\ge d$ [G1]. This is equivalent that for arbitrary
small $\lambda$ a space $X$ admits a uniformly cobounded $\lambda$-Lipschitz
map onto an $n$-dimensional uniform simplicial complex [G1].

We note that for every $n$-dimensional asymptotic polyhedron $L$,
$\as L=n$.

\proclaim{Lemma 3.1}
Suppose that $\as X\le n$ and $f:X\to\R_+$ is a given proper function.
Then there are a compact set $C\subset X$ and a 1-Lipshitz map $\phi:X\to L$ to
an $n$-dimensional asymptotic polyhedron with $diam(\phi^{-1}(\Delta))\le f(z)$
for all $z\in\phi^{-1}(\Delta)\setminus C$.
Moreover $L$ can be presented as the union
$L_0\cup J_0\cup L_1\cup J_1\cup\dots$, where each $L_i$ is a uniform
polyhedron of size $2^i$ and edges in each complex $J_i$ have the length $2^i$
or $2^{i+1}$.
\endproclaim
 \demo{Proof}
Let $\sU$ be an open cover of $X$ of multiplicity $\le n+1$ with the Lebesgue
number $> 2\lambda$ and with $\max\{\diam U\mid U\in\sU\}\le R$. We may assume
that $\max_{x\in X}d(x,X\setminus U)> 2\lambda$ for all $U\in\sU$.
Denote $\phi(x)=d(x,X\setminus U)$. The family of maps $\{\phi_U\}_{U\in\sU}$
defines an
1-Lipschitz map $\phi':X\to l_{\infty}$ to the Banach space $l_{\infty}$.
Note that $\phi'(X)\subset l_{\infty}\setminus B_{2\lambda}(0)$. Also we
note that the radial projection
$\pi$ of $l_{\infty}\setminus B_{2\lambda}(0)$ onto the sphere $S_{\lambda}(0)$
is 1-Lipschitz. Then the map $\phi=\pi\circ\phi':X\to S_{\lambda}(0)$ is 1-Lipschitz.
Since the multiplicity of the cover $\sU$ does not exceed $n+1$,
the image $\phi(X)$ lies in the $n$-skeleton of the spherical infinite dimensional
simplex $\Delta^{\infty}=S_{\lambda}(0)\cap l^+_{\infty}$ of the size $\lambda$.
In other words the nerve $N(\sU)$ of the cover $\sU$ is contained as a subcomplex
in the $n$-skeleton of $\Delta^{\infty}$. There is a constant $c_n$ which depends
on $n$ only such that the identification of an $n$-simplex in $\Delta^{\infty}$
with the standard Euclidean $n$-simplex of the size $c_n\lambda$ is
1-Lipschitz. Now we have that the map $\phi:X\to N(\sU)$ is  1-Lipschitz,
where the nerve $N(\sU)$ is given a uniform metric of size $c_n\lambda$.
 Note that $\phi^{-1}(\Delta)\subset\overline{St(x,\sU)}\subset B_R(x)$ for some
$x\in\phi^{-1}(\Delta)$  for every simplex $\Delta$.

We consider a sequence of such covers $\sU_k$ with $\lambda_k>\frac{1}{c_n}2^k$;
 $k=0,1,\dots$ and define $C=f^{-1}([0,2R_2])$. Then we construct
compact sets $C_k\supset f^{-1}([0,2R_{k+3}])$ such that $d(C_k\setminus
C_{k-1},C_{k-2})\ge R_{k+2}$. Let $\tilde\sU_k=\{U\in\sU_k\mid U\cap
(C_k\setminus C_{k-1})\ne\emptyset\}$. We form a cover $\tilde\sU=\cup_k\tilde\sU_k$.
Then the nerve $N(\tilde\sU)$ has a required type $L_0\cup J_0\cup L_1\cup J_1\cup\dots$.
It can be given the metric where each $L_i$ is a uniform complex with the size
of simplices $2^i$ and every $J_i$ is a complex with simplices of a mixed type:
their edges which formed by
two elements of the cover $\tilde\sU_i$ have the length $2^i$, and the edges
which are formed by two elements of $\tilde\sU_{i+1}$ or by one element from
$\tilde\sU_i$ and
the other element from $\tilde\sU_{i+1}$ have the length $2^{i+1}$.
Here $L_i$ is the nerve of the cover $\tilde\sU_k$ and $J_i$
is the nerve of the restriction of $\sU_k\cup\sU_{k+1}$ over the boundary
$\partial C_k$. If $x\in C_k\setminus C_{k-1}$, then $f(x)\ge 2R_{k+1}$.
Note that $\phi^{-1}(\Delta)\subset B_{R_{k+1}}(x)$ for some $x\in\phi^{-1}(\Delta)$.
Then for any $z\in\phi^{-1}(\Delta)$ we have $d(z,x)\le R_{k+1}$. Hence
$z\notin C_{k-2}$ and therefore $f(z)\ge 2R_{k+1}\ge diam(\phi^{-1}(\Delta))$.

The complex $N(\tilde\sU)$ satisfies all the requirements except it is not
necessarily $n$-dimensional.
The complexes $L_i$ are $n$-dimensional the best
estimate for dimension of complexes $J_i$ is $2n+1$.
The $n$-dimensionality can be achieved by some standard
dimension theoretic trick with the choice of covers.
We are not giving all details, since for the purpose of this paper the finite
dimensionality is enough.
\qed
\enddemo
A metric space $X$ is called {\it uniformly contractible} if there is
a function $S:\R_+\to\R_+$ such that every ball $B_r(x)$ is contractible to
a point in the ball $B_{S(r)}(x)$.

Let $x_0\in X$ be a base point. We denote $\|x\|=d_X(x,x_0)$.

\proclaim{Lemma 3.2}
Let $X$ be a uniformly contractible proper metric space with $\as X=n$.
Then given a proper function $g:X\to\R_+$ there are an $n$-dimensional
asymptotic polyhedron $N$, a proper
1-Lipschitz map $\phi:X\to N$, and a proper homotopy inverse map
$\gamma:N\to X$ with $d(x,\gamma\phi(x))<g(x)$ for all $x\in X$.
\endproclaim
\demo{Proof}
We define by induction on $i$ a lift $\gamma$ on the $i$-skeleton $N^{(i)}$
of the nerve of a cover of $X$ given by Lemma 3.1 for an appropriate choice
of $f$. We take $\gamma(v)\in\phi^{-1}(v)$ for every vertex $v$. Then using
the uniform contractibility of $X$ we can extend $\gamma$ with control over
the 1-skeleton $N^{(1)}$ and so on.
Without loss of generality we may assume that $X$ is a polyhedron of the dimension
$n$ supplied with a triangulation of mesh $\le 1$. By induction on $i$
we define a homotopy $H:X^{(i)}\times I\to X$ joining the identity map with
$\gamma\circ\phi$. We consider a function
$\psi(x)=\|x\|-\max\{d(x,y)\mid y\in H(x\times I)\}$. If $\psi$ tends to infinity,
then the map $H$ is proper. Therefore it suffices to show that $\psi$
tends to infinity for an appropriate choice of $f$. Let $S$ be a contractibility
function. We define $\rho(t)=S^{-1}(t/2)$ where $S^{-1}$ 
is the inverse function for $S$.
Then we take $f=\rho^{2n+1}\circ g$, the composition of $ g$ and $2n+1$ times
iteration of $\rho$. We assume here that $g(x)\le\|x\|/2$.
Then it is easy to verify that $\psi(x)\ge \|x\|/2$.
 \qed
\enddemo

A metric space $X$ is called {\it proper} if every ball $B_r(x)$ in $X$ is
compact.
We recall that the Higson compactification $\bar X$ of a proper metric space
$X$ is defined by the ring of bounded functions with the gradient tending to zero at
infinity [Ro1]. The reminder of this compactification is called
the {\it Higson corona} and it is denoted as $\nu X$.
Thus, $\bar X= X\cup\nu X$.
The defining property of
the Higson corona is the following:

{\it  (*) a continuous map $f:X\to Z$ to a compact
metric space is extendable over the Higson corona $\nu X$ if and only if for
every $R$ diameter of the image $f(B_R(x))$ of the $R$-ball centered at $x$
tends to zero
as $x$ approaches infinity}.

Note that a proper Lipschitz map $f:X\to Y$ induces
a continuous mapping between the Higson coronas $\bar f:\nu X\to\nu Y$.

\proclaim{Theorem 3.3}
Let $X$ be a uniformly contractible proper metric space
 with a finite asymptotic dimension
and let $\bar X$ be the Higson compactification. Then
 $\check H^{n}(\bar X;\Z_p)=0$ for all $n$ and all $p$.
\endproclaim
\demo{Proof}
We show that every map $\alpha:\bar X\to K(\Z_p,n)$ is null homotopic.
Since $\bar X$ is compact, the image $\alpha(\bar X)$ is contained in the
$k$-skeleton $K=K(\Z_p,n)^{(k)}$, $k>\as X+1$,
which is a finite complex. We fix a geodesic metric
on $K$. Let $\epsilon_K$ be an injectivity radius in $K$, i.e. every two points
within a distance $\epsilon_K$ can be joined by a unique geodesic.
 Since the map $\alpha|_X\to K$ is extendable over the Higson corona
the function $R_{\alpha}(t)=L(\alpha|_{X\setminus B_t(x_0)})$ tends to zero
at infinity. We apply Lemma 3.2 with $g(x)\le\min\{\epsilon_K/R_{\alpha}(\|x\|/2),
\|x\|/4\}$ to obtain an asymptotic polyhedron $N$ and maps $\phi:X\to N$ and
$\gamma:N\to X$. Let $[u,v]$ be an edge in $N$, then
$d_K(\alpha\gamma(u),\alpha\gamma(v))
\le R_{\alpha}(t_0)d_X(\gamma(u),\gamma(v))$ where
$t_0=\min\{\|\gamma(u)\|,\|\gamma(v)\|\}$. We may assume that there
are $x,y\in X$ such that
$\phi(x)=u$ and $\phi(y)=v$. Then $$d_X(\gamma(u),\gamma(v))\le d_X(x,y)+
\epsilon_K/R_{\alpha}(\frac{1}{2}\|x\|)+\epsilon_K/R_{\alpha}(\frac{1}{2}\|y\|)\le
diam\phi^{-1}[u,v]+2\epsilon_K/R_{\alpha}(\frac{1}{2}\|x\|)$$ provided
$\|x\|\le\|y\|$.
Because of the inequality $g(x)\le\|x\|/4$ we have that
$R_{\alpha}(t_0)\le R_{\alpha}(\|x\|/2)$. We may assume that $f<g$ and then
$diam(\phi^{-1}([u,v]))<g(x)$. Summarizing all this, we obtain the inequality
$d_K(\alpha\gamma(u),\alpha\gamma(v))\le 3\epsilon_K$. This means that
the map $\alpha\circ\gamma\circ u^{-1}$ is $3\epsilon_K$-Lipschitz where
$u:K\to K_U$ is the projection to the uniform metric.

Since $X$ is contractible, the map $\alpha\circ\gamma$ is null homotopic.
Note that the homotopy groups
$\pi_i(K)$ are finite for $i\le\dim N+1$. We apply Lemma 2.4 to obtain a
$\lambda$-Lipschitz homotopy $H:N_U\times I\to K$ of $\alpha\circ\gamma$ to
a constant map. This homotopy defines a
Lipschitz map $\tilde H:N_U\to K^I_{\lambda}$
to the space of $\lambda$-Lipschitz mappings of the unit interval $I$ to $K$.
 Then the composition
$\tilde H\circ u:N\to K^I_{\lambda}$ satisfies the Higson extendibility condition (*).
Let $\bar h:\bar N\to K^I_{\lambda}$ be the extension over the Higson corona.
This extension defines a map $\bar H:\bar N\times I\to K$. The map $\bar H$ is
a homotopy between the extension $\xi=\overline{\alpha\circ\gamma}$ and a constant map.
To complete the proof, we show that $\alpha$ is homotopic to $\xi\circ\bar\phi$
where $\bar\phi$ is the extension of the Lipschitz map $\phi$ to the Higson
compactifications. Note that
$$d_K(\alpha(x),\alpha\gamma\phi(x))\le R_{\alpha}(t_0)d(x,\gamma\phi(x))\le
R_{\alpha}(t_0)\epsilon_K/R_{\alpha}(\|x\|/2)\le\epsilon_K$$ where
$t_0=\min\{\|x\|,\|\gamma\phi(x)\|\}\ge\|x\|/2$. Then for every $x\in X$ we join
the points $\alpha(x)$ and $\alpha\gamma\phi(x)$ by the unique geodesic
$\psi_x:I\to K$.
This defines a map $\tilde\psi:X\to K^I_{\mu}$. Since both
$\alpha$ and $\alpha\circ\gamma\circ\phi$ satisfy the condition (*), the map
$\tilde\psi$ has the property (*). Let $\bar\psi:\bar X\to K^I_{\mu}$ be the
extension of $\tilde\psi$ to the Higson corona. The map $\bar\psi$ defines
a homotopy $\Psi:\bar X\times I\to K$ between $\alpha$ an $\xi\circ\bar\phi$.
\qed
\enddemo

A potential application of Theorem 3.3 to the Novikov Conjecture is
based on the following corollary. To make the connection visual, one should
compare the corollary with the Ferry-Weinberger Descent Principle
formulated in \S 5. Here (and there)  $H^{st}_*$ denotes the Steenrod homology.

\proclaim{Corollary 3.4}
Let $X$ be a uniformly contractible $n$-manifold with a finite asymptotic
dimension. Then there exists a Higson dominated metrizable corona $N$ of $X$
such that the boundary homomorphism
$H_n^{lf}(X;\Z)=\Z\to H_{n-1}^{st}(N;\Z)$ is a monomorphism.
\endproclaim
\demo{Proof}
First we apply Theorem 3.3. Then the
Schepin Spectral theorem [Dr1] implies that there is a metrizable Higson
dominated corona $N$ of $X$ such that
$\check H^*(X\cup N;\Z_p)=0$ for all primes $p$. Then for the Steenrod homology
the inclusion of the boundary induces an isomorphism
$\tilde\partial_p:H^{lf}_n(X;\Z_p)\to H_{n-1}^{st}(N;\Z_p)$ for all $p$.
Consider an $(n-2)$-connected and locally $(n-2)$-connected compactum
$Y=N\cup W$
with $\dim W\le n-1$.
Consider a diagram generated by exact sequence of pairs $(X\cup N,N)\subset
(X\cup Y,Y)$. If the image in $H^{st}_{n-1}(Y)$ of the generator
$1\in\Z=H^{lf}_n(X)$
is not divisible by $p$, then it is not divisible by $p$
in $H^{st}_{n-1}(N)$. Since $W$ is $(n-1)$-dimensional, the inclusion homomorphism
$H^{st}_{n-1}(N;\Z_p)\to H^{st}_{n-1}(Y;\Z_p)$ is a monomorphism. Hence,
the boundary homomorphism $H^{lf}_n(X;\Z_p)\to H_{n-1}^{st}(Y;\Z_p)$ is
a monomorphism.
Since $Y$ is locally $(n-2)$-connected, the Steenrod homology in the dimension
$n-1$ agrees with the singular homology and we can use the universal coefficient formula.
Since $Y$ is $(n-2)$-connected, the universal coefficient formula implies that
$H^{st}_{n-1}(Y;\Z_p)=H^{st}_{n-1}(Y)\otimes\Z_p$. Therefore
the homomorphism $H^{lf}_n(X)\otimes\Z_p\to H^{st}_{n-1}(Y)\otimes\Z_p$
is a monomorphism. Hence the image of the generator is not divisible by $p$
in $H^{st}_{n-1}(Y)$ and as well in $H^{st}_{n-1}(N)$. Therefore
$\Z=H_n^{lf}(X)\to H_{n-1}^{st}(N)$ is a monomorphism.\qed
\enddemo

REMARK.
In [Ke] Keesling established that 1-dimensional cohomology of the Higson
compactification of $\R^n$ is nonzero.
Also it was shown in [D-F] that $\check H^n(\nu\R^n;\Z)\ne 0$. Then Theorem 3.3
implies that these groups are $p$-divisible for all $p$.

\head \S 4 Hypereuclidean manifolds\endhead

Let $Y$ be a metric space with a base point $y_0$.
We define the suspension $\Sigma Y$ as the quotient metric space
of the product
$I\times Y$ with the $l_1$-metric. Here the quotient map
$q:I\times Y\to\Sigma Y$ collapses the set $\{0,1\}\times Y\cup I\times\{y_0\}$
to a point. The quotient metric is
the maximal metric with respect to which $q$ is 1-Lipschitz.
Then the
$n$-th suspension $\Sigma^nY$ of $Y$ can be identified inductively with
the following quotient map:
$$
q_n:I^n\times Y= I\times I^{n-1}\times Y\to I\times\Sigma^{n-1}Y=\Sigma^nY.
$$
 Let $f:X\to Y$ be a map. We denote by
$1_{I^n}\hat\times f=q_n\circ(1_{I^n}\times f):I^n\times X\to\Sigma^nY$.
Note that $L(1_{I^n}\hat\times f)=L(f)$.

\proclaim{Lemma 4.1 (Trading Lemma)}
For every $n$ there is a number $c_n$ such that for every $L$-Lipschitz map
$f:I^n\to Y$ which is $\lambda$-Lipschitz on the boundary $\partial I^n$,
the map $1_{I^n}\hat\times f:I^n\times I^n\to\Sigma^nY$ can be deformed
by a $c_nL$-Lipschitz homotopy $H:I^n\times I^n\times I\to\Sigma^nY$,
fixed on the boundary $\partial(I^n\times I^n)$, to a
map $g$ with $L(g|_{x\times I^n})\le c_n\lambda$  for all $x\in B^n$.
\endproclaim
\demo{Proof}
We fix a $c'_n$-Lipschitz
isotopy $H':I^n\times I^n\times I\to I^n\times I^n$ that exchanges the factors.
There is a $\bar c_n$-Lipschitz homeomorphism for some constant $\bar c_n$,
which depends on $n$ only,
$$
\xi:(I^n\times I^n)\setminus Int\frac{1}{2}(I^n\times I^n)\to
\partial(I^n\times I^n)\times I
$$
such that the restriction of $\xi$ to the boundary
$\partial(I^n\times I^n)$  identifies it with
$\partial(I^n\times I^n)\times\{0\}$
and  the restriction of $\xi$
to $\partial\frac{1}{2}(I^n\times I^n)$ is the multiplication by 2.

Let $c_n=2c_n'\bar c_n$. We define
$$
H(x,y,t)=\cases (1\hat\times f)H'((1+t)x,(1+t)y,t) & if\ \
(x,y)\in\frac{1}{2}(I^n\times I^n)\\
(1\hat\times f)(tH'(\xi(x,y))+(1-t)(x,y)) &  \ otherwise.\\
\endcases
$$
Note that $$L((1\hat\times f)H'((1+t)x,(1+t)y,t))\le 2c_n'L\le c_nL$$ and
$$L((1\hat\times f)(tH'(\xi(x,y))+(1-t)(x,y)))\le L(c_n'\bar c_n+1)\le c_nL.$$
Thus, the map $H$ is $c_nL$-Lipschitz.

It is easy to check that $H(x,y,0)=1\hat\times f(x,y)$ for all $x,y\in I^n$.

If $(x,y)\in\partial(I^n\times I^n)$
we have that
$$
H(x,y,t)=(1\hat\times f)(tH'(\xi(x,y))+(1-t)(x,y))=
(1\hat\times f)(t(x,y)-(1-t)(x,y))=(1\hat\times f)(x,y).
$$
Thus, the homotopy
$h_t$ is fixed on $\partial(I^n\times I^n)$.

Consider the map
$$
g(x,y)=H(x,y,1)=\cases (1\hat\times f)(2y,2x) & \ if\ \
(x,y)\in\frac{1}{2}(I^n\times I^n)\\
(1\hat\times f)H'(\xi(x,y)) &\ \ otherwise.\\
\endcases
$$
The Lipschitz constant of $(1\hat\times f)(2y,2x)$ as a function of $y$ is equal
2. Since $L((1\hat\times f)|_{\partial(I^n\times I^n)})\le\lambda$, we have
that $L((1\hat\times f)\circ H'\circ\xi)\le c_n'\bar c_n\lambda\le c_n\lambda$.
Hence $L(g|_{x\times I^n})\le c_n\lambda$.
\qed
\enddemo
For a map $f:X\times Z\to Y$ between metric spaces we denote
$$
L^Z(f)=\sup\{L(f|_{X\times z})\mid z\in Z\}, \ \ \ \
L^X(f)=\sup\{L(f|_{x\times Z})\mid x\in X\}.
$$

The loop space $\Omega Y$ on a metric space $Y$ as any other mapping space
is endowed with the natural $\sup$ metric. This makes the iterated loop
space $\Omega^nY$ equal to the mapping space $(Y,y_0)^{(I^n,\partial I^n)}$.
Every map $f:I^n\times Z\to Y$ with $f(\partial I^n\times Z)=y_0$ induces
a map $F:Z\to\Omega^n Y$ by the formula $F(z)(x)=f(x,z)$. If
$L^{I^n}(f)\le\lambda$, then $L(F)\le\lambda$.

For every $n$ and every $Y$ there is a natural inclusion $j_n^Y:Y\to\Omega^n\Sigma^nY$
defined by the following rule: A point $y\in Y$ corresponds under $j_n^Y$ to
the map $F_y=q_n\circ(1_{I^n}\times c_y):I^n\to\Sigma^nY$ where $c_y$ is
a constant map $c_y:I^n\to Y$ to $y$.
One can check that $j_n^Y$ is an isometric imbedding.
Generally, one can define a map
$j_{i,n+i}^Y:\Omega^i\Sigma^iY\to\Omega^{n+i}\Sigma^{n+i}Y$ by taking a map
$\phi:I^n\to\Sigma^nY$ to the map $q_i(1_{I^i}\times\phi)$.
Then $j_n^Y=j_{0,n}^Y$.
Similarly $j_{i,n+i}^Y$ is an isometry.
We note that
$j^Y_{i,n+k+i}=j^Y_{k+i,n+k+i}\circ j^Y_{i,k+i}$ for all $i$, $k$, and $n$.

\proclaim{Corollary 4.2}
For every $L$-Lipschitz map $f:I^n\to Y$ with a $\lambda$-Lipschitz restriction
on the boundary $\partial I^n$, the composition
$j_n^Y\circ f:I^n\to\Omega^n\Sigma^nY$ is homotopic to
a $c_n\lambda$-Lipschitz map by means of $c_nL$-Lipschitz homotopy $H_t$
which is fixed on $\partial I^n$.
\endproclaim
\demo{Proof}
Consider a homotopy $h_t:I^n\times I^n\to\Sigma^nY$ from Lemma 4.1.
It defines a $c_nL$-Lipschitz homotopy $H_t:I^n\to\Omega^n\Sigma^nI^n$
by the formula $H_t(x)(z)=h_t(x,z)$. Then
$$
H_0(x)(z)=h_0(x,z)=(i\hat\times f)(x,z)=q_n(x,f(z))=F_{f(z)}(x)=j_n^Y(f(z))
$$
Since $H_1(x)(z)=g(x,z)$ and $L^{I^n}(g)\le c_n\lambda$, we have
$L(H_1(x))\le c_n\lambda$. Since the homotopy $h_t$ keeps the set
$\partial(I^n\times I^n)$ fixed, the homotopy $H_t$ keeps the boundary
$\partial I^n$ fixed.\qed
\enddemo

By $\Omega^n_LW$ we denote the subset of $\Omega^nW$ consisting of
$L$-Lipschitz maps $\phi:I^n\to W$.

For any numbers $n$ and $i$ denote by $\alpha_{n,i}^Y$ a map from
$\Omega^n\Sigma^n\Omega^i\Sigma^iY \to \Omega^{n+i}\Sigma^{n+i}Y$
defined by the following formula:
$\alpha_{n,i}^Y(\Phi)(z,x)=q_{n+i}((y(z),x),\phi_z(x))$ where
$\Phi(z)=q_n(y(z),\phi_z)$, $z\in I^n$, $x\in I^i$, $y:I^n\to I^n$, and
$\phi_z:I^i\to\Sigma^iY$.

\proclaim{Proposition 4.3}
The map $\alpha_{n,i}^Y$ is 1-Lipschitz.
\endproclaim
\demo{Proof}
By the definition we have

$dist(\alpha_{n,i}^Y(\Phi),\alpha_{n,i}^Y(\Phi'))=\sup_{x,z}
dist(\alpha_{n,i}^Y(\Phi)(z,x),\alpha_{n,i}^Y(\Phi')(z,x))\le$

$\sup_{x,z}dist((y(z),x),\phi_z(x)),(y'(z),x),\phi'_z(x))=
\sup_{x,z}\|y(z)-y'(z)\|+\sup_{x,z}\|\phi_z(x)-\phi_z'(x)\|\le
dist(\Phi,\Phi')$.
\qed
\enddemo
\proclaim{Proposition 4.4}
$j^Y_{i,n+i}=\alpha^Y_{n,i}\circ j_n^{\Omega^i\Sigma^iY}$
for all $i$, $n$ and $Y$.
\endproclaim
\demo{Proof}
Let $\phi\in\Omega^i\Sigma^iY$. Then $j^Y_{i,n+i}(\phi)=q_i(1_{I^n}\times\phi)$.
On the other hand $
\alpha^Y_{n,i}\circ j_n^{\Omega^i\Sigma^iY}(\phi)
=\alpha^Y_{n,i}(q_n(1_{I^n},c_{\phi}))=q_{n+i}(1_{I^n}\times
1_{I^i},\times\phi)=q_i(1_{I^n}\times\phi)$. \qed
\enddemo

\proclaim{Lemma 4.5}
For every natural number $n$ there exists a number $b_n$ such that
for any positive $\lambda$, any two tending to infinity functions
$\psi:\R_+\to\R_+$, $\xi:\R_+\to\R_+$ and every
continuous map $f:K\to\Omega^i\Sigma^iY$ of
an $n$-dimensional uniform polyhedron with $L_f(t)\le\psi(t)$,
$L(f|_{K^{(n-1)}})\le \lambda$,
and $f(x)\in\Omega^i_{\xi(\|x\|)}\Sigma^iY$ for all $x\in K$,
there exists a homotopy
$H:K\times I\to\Omega^{n+i}\Sigma^{n+i}Y$ deforming a map $j_{i,n+i}^Y\circ
f$ to a map $g$ with the properties:  \roster \item{} $L(g)\le b_n\lambda$;
\item{} $L_H(t)\le b_n\psi(t)$ and
\item{} $g(x)\in\Omega^{n+i}_{2b_n(\psi(\|x\|)+\xi(\|x\|))}\Sigma^{n+i}Y$.
\endroster
\endproclaim
\demo{Proof}
We fix an $n$-simplex $\Delta\subset K$ and consider the restriction
$f_{\Delta}$ of $f$ to it. Note that
 $L(f_{\Delta}|_{\partial\Delta})\le\lambda$
and $L(f_{\Delta})\le\psi(\|\Delta\|)$, where
$\|\Delta\|=\max\{dist(x,x_0)\mid x\in\Delta\}$.
Since $\Delta$ is $a_n$-Lipschitz homeomorphic to the cube $I^n$, we may apply
Corollary 4.2. According to Corollary 4.2 there is a
$a_nc_nL(f_{\Delta})$-Lipschitz
homotopy in $\Omega^n\Sigma^n(\Omega^i\Sigma^iY)$ from
$j_n^{\Omega^i\Sigma^iY}\circ f$ to
a $a_nc_n\lambda$-Lipschitz map $g_{\Delta}$ which is fixed on the boundary
$\partial\Delta$. The union of these homotopies defines a homotopy
$\tilde H:K\times I\to\Omega^n\Sigma^n(\Omega^i\Sigma^iY)$ from
$j_n^{\Omega^i\Sigma^iY}\circ f$ to a $a_nc_n\lambda$-Lipschitz map
$\tilde g$ satisfying
the inequality $L_{\tilde H}(t)\le a_nc_n\psi(t)$.
We define $H=\alpha^Y_{n,i}\circ\tilde H$.
Then $H$ is a homotopy between $g=\alpha^Y_{n,i}\circ\tilde g$ and
$\alpha^Y_{n,i}\circ j_n^{\Omega^i\Sigma^iY}\circ f$.
The latter is equal to $j^Y_{i,n+i}\circ f$ by Proposition 4.4.
The condition (1) holds for $\tilde g$ with $b_n=a_nc_n$ by
the construction. Then by Proposition 4.3 it holds for $g$.
Similarly, the condition (2)  holds for $\tilde H$ and by Proposition 4.3
it holds for $H$.

Let $\tilde g(x)=\Phi\in\Omega^n\Sigma^n\Omega^i\Sigma^iY$. Note that
$\Phi(z)=q_n(y(z),\phi_z)$,
where $y:I^n\to I^n$ and $\phi_z:I^i\to\Sigma^iY$ are two maps.
By the construction of $\tilde g$ we have $\tilde g(x)\in
\Omega_{b_n\psi(\|x\|)}^n\Sigma^n\Omega^i_{\xi(\|x\|)}\Sigma^iY$.
It means that the maps $y$ and $\phi_z$ are $b_n\psi(\|x\|)$-Lipschitz
as functions of $z$. Additionally, the map $\phi_z$ is $\xi(\|x\|)$-Lipschitz
for every $z$.
Let $(z,t), (z',t')\in I^n\times I^i=I^{n+i}$.
Then

$dist(\alpha^Y_{n,i}(\Phi)(z,t),
\alpha^Y_{n,i}(\Phi)(z',t'))\le
dist(q_{n+i}(y(z),t),\phi_z(t)),q_{n+i}(y(z'),t'),\phi_{z'}(t'))$

\

$\le
dist((y(z),t),\phi_z(t)),(y(z'),t'),\phi_{z'}(t'))=\|y(z)-y(z')\|+
\|t-t'\|+dist(\phi_z(t),\phi_{z'}(t'))$

\

$\le b_n\psi(\|x\|)\|z-z'\|+
\|t-t'\|+dist(\phi_z(t),\phi_z(t'))+dist(\phi_z(t'),\phi_{z'}(t'))\le$

\

$\max\{b_n\psi(\|x\|),1\}\|(z,t)-(z',t')\|+\xi(\|x\|)\|t-t'\|+
b_n\psi(\|x\|)\|z-z'\|\le $

\

$2b_n(\psi(\|x\|)+\xi(\|x\|)).$
Here we assume that $\psi\ge 1$.
Thus, the condition (3) is checked.
\qed
\enddemo

Let $f:K\to Y$ be a map of an infinite simplicial complex to a metric space
which is Lipschitz on every finite complex and let $x_0\in K$ be a base point.
We define
 $L_f((t)=L(f|_{B_t(x_0)})$, where $B_t(x_0)$ is
a ball
of radius $t$ centered at $x_0$. For a homotopy $H:K\times I\to Y$ 
we denote by $L_H(t)=L(f|_{B_t(x_0)\times I})$.

Let $\phi:\R_+\to\R_+$ be a function, by $B^n_{\phi}\times K$ we denote
the set $\{(z,x)\in\R^n\times K\mid \|z\|\le\phi(d_K(x,x_0))\}$.

\proclaim{Lemma 4.6}
Let $K$ be a finite dimensional uniform polyhedron and let
$f:K\to Y$ be a continuous map to a compact polyhedron $Y$ with
$L_f(t)\le\psi(t)$ for some function $\psi$.
Then there are
numbers $n,a$ and $\lambda$
and a map $g: K\to \Omega^n\Sigma^nY$ with
$g(x)\in\Omega_{a\psi(\|x\|)}^n\Sigma^nY$ satisfying the inequalities
$L(g)<\lambda$ which is homotopic to the map $j^Y_n\circ f$ by means of a
homotopy $H$ with the property $L_H(t)\le a\psi(t)$.  \endproclaim
\demo{Proof}
We apply induction on $\dim K=i$.
If $\dim K=0$, then we assume that all points in $K$ are at least
one apart. Then the statement follows from compactness of $Y$ with
$n=0$, $a=1$ and $\lambda=2diam(Y)$.

Assume that the lemma holds for $i$ and let $\dim K=i+1$.
Let $f:K\to Y$ be given with $L_f(t)<\psi(t)$.
By the inductive assumptions, there is a map
$g_i:K^{(i)}\to\Omega^r\Sigma^rY$ with $L(g_i)<\lambda_i$,
$g_i(x)\in\Omega_{a_i\psi(\|x\|)}^r\Sigma^rY $, and a homotopy $H^i$, joining
$j^Y_r\circ f$ with $g_i$, which satisfies the inequality $L_{H^i}(t)\le a_i\psi(t)$. We
consider an arbitrary $(i+1)$-simplex $\Delta\subset K$ and consider the map
$f$ restricted to $\Delta$. Then the homotopy $H^i$ on
$\partial\Delta$ can be extended to a homotopy $\bar H^i$ of
$\Delta$ by means of an isotopy $G:\Delta\times I\to
\Delta\times\{0\}\cup \partial\Delta\times I$,  which joins
the embedding $\Delta$ as $\Delta\times\{0\}$ to a homemomorphism between
$\Delta$ and $\Delta\times\{0\}\cup\partial\Delta\times I$. We may assume
that $G$ is $c$-Lipschitz where $c$ depends on $i$ only.
Let $f_1=\bar H^i_1:\Delta\to\Omega^r\Sigma^rY$. Note that
$L(\bar H^i)<cL_f(d(\Delta,x_0))$ and
$L(f_1|_{\partial\Delta})=L(g_i)<\lambda_i$.
Let $\tilde H$ be the union of homotopies $\bar H^i$ for all
$i+1$-simplices $\Delta$ and let $\tilde f=\tilde H_1$.
We note that $L_{\tilde H}(t)\le cL_f(t)\le c\psi(t)$ and $L(\tilde
f|_{K^{(i)}})<\lambda_i$.
Also we note that $\tilde f(x)\in\Omega_{a_i\psi(\|x\|)}^r\Sigma^rY$.
We apply Lemma 4.5 to the map $\tilde f:K\to\Omega^r\Sigma^rY$ with
$n=i+1$ to obtain a homotopy $H'$ transforming $j^Y_{r,i+1+r}\tilde f$ to
a map $g:K\to\Omega^{i+1+r}\Sigma^{i+1+r}Y$ with the properties:
(1) $L(g)\le b_n\lambda_i$; (2) $L_{H'}(t)\le b_nc\psi(t)$;
(3) $g(x)\in\Omega^{i+1+r}_{2b_n(1+a_i)\psi(\|x\|)}\Sigma^{i+1+r}Y$.
 We set $\lambda=b_n\lambda_i$ and $a=2b_n\max\{2(1+a_i),c\}$ and
define a homotopy $H$ as $j^Y_{r,i+1+r}\tilde H$ followed by $H'$.
Note that $H$ joins $j^Y_{r,i+1+r}\circ j^Y_r\circ f=j^Y_{i+1+r}$ with $g$.
Then all conditions are satisfied.
\qed
\enddemo
\proclaim{Corollary 4.7}
Let $K$ and $f:K\to Y$ be as above with  $L_f(t)\psi(t)$
where $\psi$ is a Lipschitz function. Then there
are a number $n$ and a Lipschitz map
$q:B^n_{\psi}\times K\to\Sigma^nY$ such that $q$ is homotopic to
$(1\hat\times f)\circ(\frac{1}{\psi(\|\ \|)}\times 1_K)$ by means of a homotopy $h_t$
with $h_t(\partial B^n_{\psi}\times K)=y_0$.
\endproclaim
\demo{Proof} Let $g:K\to\Omega^n\Sigma^nY$
 be given by Lemma 4.6. The map $g$ defines a map $\hat g:B^n\times
K\to\Sigma^nY$ with the properties:
(1) $L^{B^n}(\hat g)\le\lambda$; (2) $L(\hat g|_{B^n\times\{x\}})\le
 a\psi(\|x\|)$; (3) $\hat g$ is homotopic to $1\hat\times f$ by means of
 homotopy $h'_t$ with $h'_t(\partial B^n_{\psi(t)}\times K)=y_0$.

 We define $q(b,x)=\hat g(\frac{b}{\psi(\|x\|)},x)$, where $x\in K$, 
$b\in B^n_{\psi(\|x\|)}$.  We show that
$q$ is Lipschitz.  First we consider a pair of points $(b_1,x),(b_2,x)\in
B^n_{\psi}\times K\subset\R^n\times K$. Note that

$$d(q(b_1,x),q(b_2,x))\le d(\hat g(b_1/\psi(\|x\|),x),\hat 
g(b_2/\psi(\|x\|),x))\le L_g(\|x\|)\|b_1/\psi(\|x\|)-b_2/\psi(\|x\|)\|.$$

Then by the conclusion of Lemma 4.6 we have

\

$d(q(b_1,x),q(b_2,x))\le a\psi(\|x\|)\|\frac{b_1}{\psi(\|x\|)}-
\frac{b_2}{\psi(\|x\|)}\|\le a\|b_1-b_2\|$.

\

Next, we consider a pair of type $(b,x_1)$ and $(b,x_2)$. Then
$$d(q(b,x_1),q(b,x_2))=d(\hat g(\frac{b}{\psi(\|x_1\|)},x_1),\hat g(
\frac{b}{\psi(\|x_2\|)},
x_2))\le$$
$$d(\hat g(\frac{b}{\psi(\|x_1\|)},x_1),\hat g(\frac{b}{\psi(\|x_1\|)},
x_2))+d(\hat g(\frac{b}{\psi(\|x_1\|)},x_2),\hat g(
\frac{b}{\psi(\|x_2\|)},x_2))$$
$$\le \lambda d(x_1,x_2)+a\psi(\|x_2\|)\|b/\psi(\|x_1\|)-b/
\psi(\|x_2\|)$$
$$\le \lambda d(x_1,x_2)+a|\psi(\|x_2\|)-\psi(\|x_1\|)|
\|\le (\lambda+a\mu)d(y_1,y_2)$$
where $\mu$ is
 a Lipschitz constant for $\psi$.
\qed
\enddemo

Let $(K,d)$ be an asymptotic polyhedron and let $u_K:K\to K_U$ denote
the identity map which switch from the metric $d$ to the uniform metric.
Let $f:\R_+\to\R_+$ be a tending to infinity function, we say that $K$
has growth less than $f$ if $\|x\|_d <f(\|x\|_U)$. We say that $K$
has a {\it quadratic growth} if $\|x\|_d<\|x\|_U^2$.

\proclaim{Lemma 4.8}
Let $X$ be the universal cover of an aspherical manifold $M_0$ ,  $p:X\to M_0$,
supplied with the induced metric,
and let $\as X=n$.
Then, given a tending to infinity function $\beta:\R_+\to\R_+$, there are
an $n$-dimensional asymptotic polyhedron $N$, a proper
1-Lipschitz map $\phi:X\to N$, and a proper homotopy inverse map
$\gamma:N\to X$ with $L_{\gamma\circ u^{-1}_K}(t)<\beta(t)$, where $K\subset N$
and the closure of $N\setminus K$ is compact.
\endproclaim
\demo{Proof}
 First we show that for the projection $\phi_{\sU}:X\to N(\sU)$ to the
uniform nerve of a uniformly bounded cover $\sU$ admits a Lipschitz
homotopy inverse $\gamma:N(\sU)\to X$. We define $\gamma$ on the vertices
$v$ by taking points $\gamma(v)\in\phi_{\sU}^{-1}(v)$ and define $\gamma$
on every edge $[v,u]$ by sending it with a constant speed along a minimal
geodesic joining
$\gamma(v)$ and $\gamma(u)$. Assume that a $\lambda_i$-Lipschitz map
$\gamma:N^{(i)}(\sU)\to X$ is already constructed for $i\ge 1$.
For every $(i+1)$-dimensional simplex $\Delta $ we consider the map
$f=p\circ\gamma|_{\Delta}$. Since this map is null homotopic and
$\pi_{i+1}(M_0)=0$, Lemma 2.3 implies that there exists
$\lambda_{i+1}$-Lipschitz extension $\xi_{\Delta}$, $\lambda_{i+1}=\mu(\lambda_i)$.
Then there is a unique
lift $\xi'$ of $\xi$ that extends $\gamma|_{\partial\Delta}$. Since $p$ is a
local isometry, the lift $\xi'$ is also $\lambda_{i+1}$-Lipschitz. The union
of all $\xi'$ defines a $\lambda_{i+1}$-Lipschitz map on the $(i+1)$-skeleton.

 For a homotopy lift $\gamma: N\to X$ from Lemma 3.2 we define a function
$\psi_{\gamma}(t)=L(\gamma|_{B_t(x_0)})$. Because of the above remark we
always may assume that the function $\psi_{\gamma}$ takes finite values.
Moreover, for any given tending to infinity function $\beta:X\to\R_+$
there are an approximation $\phi:X\to N$ by an asymptotic
$n$-dimensional polyhedron and a proper homotopy inverse map $\gamma:N\to X$
with $\psi_{\gamma}<\beta$.

 For every approximation by an asymptotic polyhedron one can find better
approximation by an asymptotic polyhedron of the quadratic growth.
To do that we enlarge the sets $C_i$ in the proof of Lemma 3.1 so that
$C_i$ contains the set $\phi_i^{-1}(B_{2^{i}}(x_0))$, where
$\phi_i:X\to N(\sU_i)$ is the projection to the nerve with the uniform metric.
Then for every $x\in B_{2^{i+1}}(x_0)\setminus B_{2^{i}}(x_0)$ we have
$\|x\|_U=\Sigma_{k=1}^i2^{i}+c$ where $c$ is the length of a segment in
$ B_{2^{i+1}}(x_0)\setminus B_{2^{i}}(x_0)$.
We note that $$\|x\|_d\le\Sigma_{k=1}^i2^i2^i+2^{i+1}c\le
 (\Sigma_{k=1}^i2^i+c)^2=\|x\|_U^2.$$
In view of this we can find an arbitrarily close approximation $\phi:X\to N$ by
asymptotic polyhedra of the quadratic growth. According to the above remark 
we can take
an approximation such that
there exists a homotopy inverse map
$\gamma:N\to X$ with $L_{\gamma}(t)<\beta(\sqrt{t})$. Then
$$L_{\gamma\circ u^{-1}_K}(t)=
L(\gamma\circ u^{-1}_K|_{B^U_t(x_0)})\le L(\gamma|_{u^{-1}_K(B^U_t(x_0))})\le
L(\gamma|_{B_{t^2}(x_0)})=L(\gamma)(t^2)<
\beta(t)$$ \qed
\enddemo

DEFINITION. An $m$-manifold $M$ is called integrally {\it hypereuclidean}
if there is a proper Lipschitz map $f:M\to\R^m$ of degree one.
\proclaim{Theorem 4.9}
Let $M$ be the universal cover of a closed aspherical $m$-manifold with 
$\as M<\infty$.
Then the product $M\times R^n$ is integrally hypereuclidean for some
$k$.
\endproclaim
\demo{Proof}
There is a Lipschitz map $\alpha:(M\times\R)\setminus C\to S^m$ of 
degree one for some compact set $C$. Indeed, we may assume that 
for some ball $B_r$ in $M$ there is a Lipschitz map
$\gamma:(B_r,\partial B_r)\to (B_r/\partial B_r,\{\partial B_r\})\to (S^m,y_0)$ 
to
the unit sphere which is a relative homeomorphism. Let $C=B\times[-1,1]$. 
One can take
$$
\alpha(x,t)=\cases \gamma(x) &\ \ if\ \ x\in B_r\ \ and\ \ t\ge 1\\
                    y_0 &\ \ otherwise\\
\endcases
$$

Let $L(\alpha)<a$. Let $\psi:\R_+\to\R_+$ be a Lipschitz map of a sublinear
growth.

Let $\phi:M\times\R\to N$ be an approximation by an asymptotic polyhedron
of quadratic growth
as in Lemma 4.8 with a homotopy inverse map $\gamma:N\to M\times\R$ with
$L_{\gamma}(t)\le\max\{\psi(\sqrt{t})/a,c\}$.
We can find a subpolyhedron
$K\subset L$ with a relatively compact complement such that
$L_{\gamma|_K}(t)\le\psi(\sqrt t)/a$ and the map
$\alpha\circ\gamma|_K\circ\phi|_{\gamma^{-1}(K)}$
induces the same element $[\alpha|_{\gamma^{-1}(K)}]$ in the $m$-dimensional
cohomology. Let $f=\alpha\circ\gamma|_K\to S^m$,
where $K$ is considered as a uniform complex.  Then $L_f(t)\le \psi(t)$.
Then by the Corollary 4.7 there is $n$ and a Lipschitz map
$q:B^n_{\psi(t)}\times K\to S^{n+m}$ such that  $q$ is homotopic to
$(1\hat\times f)\circ(\frac{1}{\|\ \|}\times 1_K)$ by means of a homotopy $h_t$
with $h_t(\partial B^n_t\times K)=y_0$.

We define a map $\xi:\R^n\to\R^n$ as $\xi(z)=z/\psi(\|z\|)$ and
 consider the map $g=q\circ(\xi|_{B^n_t}\times u_K):B^n_t\times K\to S^{n+m}$
where $K$ is taken with its metric of  an asymptotic polyhedron.
It is easy to see that the map $g$ has the Lipschitz constant vanishing at
infinity. We extend the map $g$ to a map $\bar g:\R^n\times K\to S^{n+m}$
by the constant map. Then the map
$$w=\bar g\circ(1_{\R^n}\times\phi|):\R^n\times
((M\times\R)\setminus C)\to S^{n+m}$$
can be continuously extended to $\bar w$ by the constant map over
$(\R^n\setminus B_R^n)\times C$ for a large enough $R$. Then $\bar w$
has the gradient tending to zero at infinity. Therefore $w$
can be extended over the Higson corona $\nu(\R^{n+1}\times M)$. We note that
the map $\bar w$ has degree one. Then it follows that $\R^k\times M$ is
integrally hypereuclidean for $k=n+1$ [Ro1].
\qed
\enddemo

\head \S 5 Application to the Novikov Conjecture \endhead

There are several different coarse type conditions on a manifold $M$ that
apply the Novikov Higher Signature Conjecture. A comparison of some of them
was attempted in [Dr1]. We will use here the following refinement of
Carlson-Pedersen's conditions [C-P], [F-W] which is due to
Ferry-Weinberger:
\proclaim{ Descent Principle [D-F]} Let $p:X\to M$ be the
universal cover of a closed aspherical $n$-manifold with the fundamental
group $\Gamma=\pi_1(M)$. Suppose that $X$ admits a $\Gamma$-equivariant
metrizable compactification
$\bar X= X\cup\partial X$ such that the boundary homomorphism
$H^{lf}_n(X;\Q)\to H^{st}_{n-1}(\partial X;\Q)$ is a $\Gamma$-equivariant
split monomorphism for the Steenrod homology. Then the Novikov conjecture
holds for
$M$.
\endproclaim
As we mentioned in Introduction,
G. Yu proved the coarse Baum-Connes conjecture for all geometrically finite
groups with finite asymptotic dimension [Yu1]. In particular his result implies
the Novikov conjecture for these groups [Ro2]. Here we give an alternative
prove of the Novikov conjecture based on the Ferry-Weinberger Descent Principle.

\proclaim{Theorem 5.1}
Suppose that the fundamental group $\Gamma=\pi_1(M)$ of a closed
aspherical manifold $M$ has a finite asymptotic dimension as a
metric space with the word metric. Then the Novikov conjecture holds for $M$.
\endproclaim

We note that an action of a group $\Gamma$ by isometries on a metric space $X$
can be extended over the Higson corona.

DEFINITION. An $n$-manifold $X$ with an isometric action of a group $\Gamma$
on it is called $\Gamma$-{\it hypereuclidean} if its Higson corona $\nu X$
admits
a map $f:\nu X\to S^{n-1}$ of degree one which induces
a $\Gamma$-equivariant homomorphism of the cohomologies for the trivial action
 on $S^{n-1}$.
In other words $f$ defines a cohomology element fixed under the action of $\Gamma$
on $\check H^{(n-1)}(\nu X)$.
By the definition the degree of $f$ is the degree
of the homomorphism

$$\Z\to H^{n-1}(S^{n-1})@>f^*>>\check H^{n-1}(\nu X)@>\delta>> H^n_c(X)=\Z.$$

The proof of Theorem 5.1 is based on the following two lemmas.
\proclaim{Lemma 5.2} Suppose that the universal cover $X$ of a closed aspherical
$m$-manifold $M$ is equivariantly hypereuclidean. Then the Novikov conjecture
holds for $M$.
\endproclaim
\demo{Proof}
Let $f:X\to\R^m$ be a proper Lipschitz map of degree one such that
the extension map $\bar f$ to the Higson corona 
is a cohomology equivariant map to $S^{m-1}$.
Consider the $\Gamma$-equivariant commutative diagram induced
by this extension and the boundary maps
$$
\CD
\check H^{m-1}(\nu X;\Q) @>\delta>> H^m_c(X;\Q)\\
@A{\bar f^*}AA  @Af^*AA\\
H^{m-1}(S^{m-1};\Q) @>\tilde\delta>> H^m_c(\R^m;\Q)\\
\endCD
$$
Here the homomorphism $\delta\circ\bar f^*$ is an equivariant isomorphism.
Therefore the homomorphism $\delta$ is a $\Gamma$-equivariant split epimorphism.

Then by the Shchepin spectral theorem  one can find a
$\Gamma$-equivariant metrizable compactification $\bar Y=Y\cup\partial Y$
such that the boundary homomorphism $\bar\delta:\check H^{m-1}(\partial Y;\Q)\to
H^m_c(Y;\Q)$ is a $\Gamma$-equivariant split monomorphism
(see the proof of Lemma 9.3 in [Dr1]).  Then the boundary
homomorphism for the Steenrod homology
$$\bar\partial:H^{lf}_m(Y;\Q)\to H^{st}_{m-1}(\partial Y;\Q)$$
is a $\Gamma$-equivariant split monomorphism.
The Descent Principle completes the proof.\qed
\enddemo
\proclaim{Lemma 5.3}
Let $\Gamma=\pi_1(M)$ be the fundamental group of an aspherical $m$-manifold
and let $X$ be its universal cover. Assume that $\as X<\infty$.
Then there is a number $N$ such that $X\times\R^N$ is $\Gamma$-hypereuclidean
\endproclaim
\demo{Proof}
 We fix an isometry $h:X\to X$ that
preserves the orientation. We consider a Lipschitz map
$\alpha:X\times\R\setminus C\to S^m$ from the proof of Theorem 4.9.
First, we note that $\alpha|_{X\times\R\setminus C'}$ is Lipschitz homotopic to
$\alpha\circ (h\times 1)|_{X\times\R\setminus C'}$, where $C'$ is a compact set
in $X\times\R$ that contains $C$ and $h^{-1}(C)$.
The proof of Theorem 4.9 produces a map of the degree one
$w:\R^n\times(X\times\R\setminus C')\to S^{n+m}$ extendable over the 
Higson corona, and a Lipschitz homotopy $H$, connecting $1\hat\times\alpha$ 
with $w$.
Then $H\circ(1_{\R^n}\times h\times 1_{\R})$ is a Lipschitz homotopy 
connecting
$1\hat\times\alpha\circ(h\times 1)$ with $w\circ(1_{\R^n}\times 
h\times 1_{\R})$. It means that the maps $w$ and $w\circ(1_{\R^n}\times h\times 
1_{\R})$
are Lipschitz homotopic. Then, like in the proof of Theorem 4.9, we can show
that there is a Lipschitz homotopy $H'$ between $1_{\R^k}\hat\times w$ and
$1_{\R^k}\hat\times w\circ(1_{\R^n}\times h\times 1_{\R})$ for some $k$
which is extendable over the Higson corona for every $t\in[0,1]$ .
This homotopy defines a map 
$$G:\R^{k+n}\times(X\times\R\setminus C')\to
(S^{k+n+m})^I_{\lambda}$$ 
to the space of all $\lambda$-Lipschitz mappings of the
interval $I$ to the sphere $S^{k+n+m}$ for some $\lambda$
which is compact. Then it follows 
that
$G$ is extendible over the Higson corona. Therefore, $G$ defines a homotopy
$\bar G:\nu(X\times\R^{n+k+1})\times I\to S^{k+n+m}$ between
the extensions $\tilde w$ and $\tilde w_h$ of the maps $1\hat\times w$ 
and $1\hat\times(w\circ(1\times h)$
over the Higson corona. This means that the map $\tilde w$ is 
cohomologically $h$-invariant.

Let $h_i:X\to X$, $i=1,\dots, l$, be isometries on $X$ defined by generators 
of $\Gamma$. We apply the above argument to each  $h_i$ and take maximal $k$.
Then $N=k+n+1$.
\qed 
\enddemo
\demo{Proof of Therem 5.1}
Since the action of $\Z^N$ on $X\times\R^N$ induces the trivial action 
on $\nu(X\times\R^N)$, the product $X\times\R^N$ is $\Gamma'$-hypereuclidean 
for $\Gamma'=\Gamma\times\Z^N$. Then by Lemma 5.2 the Novikov
Conjecture holds for $\Gamma'$ and hence for $\Gamma$.\qed
\enddemo 

\proclaim{Corollary 5.4}
Suppose that the classifying space $B\Gamma$ is a finite complex 
for a group $\Gamma$ with $\as\Gamma<\infty$. Then the Novikov conjecture holds for 
$\Gamma$.
\endproclaim
\demo{Proof}
For every group $\Gamma$ with finite $B\Gamma$ M. Davis [D] gave a 
construction of an aspherical closed manifold $M$ such that $B\Gamma$ is
a retract of $M$. In [B-D] (Theorem 8) we proved that
$\as M=\as\pi_1(M)<\infty$ provided $\as\Gamma<\infty$.

It is well-known that the Novikov Conjecture for $\pi_1(K)$ 
is equivalent to 
the injectivity of the assembly map $A_K$ [F-R-R]. Hence, by Theorem 5.1,
$A_M$ is a monomorphism.
We consider the
following diagram:
$$
\CD
H_*(B\pi_1(K);\Q) @>A_K>> L_*(\Z[\pi_1(K)])\otimes\Q\\
@ Vi_*VV   @VVV\\
H_*(B\pi_1(M);\Q) @>A_M>> L_*(\Z[\pi_1(M)])\otimes\Q\\
\endCD
$$
 Since the inclusion $i:K\to M$ 
admits
a retraction, it follows that $i^*$ is a monomorphism. 
Since $A_K$ is a left 
divisor of a monomorphism $A_M\circ i_*$, it is a monomorphism itself.
\qed
\enddemo
\

In the conclusion of this section we formulate a problem which can be 
considered as a stable
version of the Higson conjecture [Ro1]. We recall that the original 
unstable version is incorrect even in the cases where the Novikov 
conjecture is true
[Ke].

Let $$p:\overline{X\times\R}=X\times\R\cup\nu(X\times\R)\to 
X\times\R\cup\Sigma\nu X=\Sigma\bar X$$ be the natural projection 
of compactifications of $X\times\R$.
We define a group 

$$SH^i(X)=\lim_{\to}\{\check H^i(\bar X) @>\Sigma>>\check H^{i+1}(\Sigma\bar X) 
@>p^*>>\check H^i(\overline{X\times\R})\to\dots\}.$$

\proclaim{Problem 5.5}
Is $SH^*(X)=0$ for the universal cover $X$ of a closed aspherical manifold $M$?
\endproclaim

We note that the affirmative answer to this problem implies the Novikov conjecture
for $M$.

\head \S 6 Asymptotically piecewise Euclidean spaces and the Gromov-Lawson 
conjecture \endhead

We recall that a metric space $X$ has a {\it bounded geometry} if for every
$\epsilon>0$ and $R>0$ there is a number $C$ such that for every $x\in X$
every $\epsilon$-net in the ball $B_R(x)$ contains no more than $C$ points.
A finitely generated group with the word metric or the universal cover
of a finite complex with the lifted metric are typical examples of metric spaces
of bounded geometry. 
Every metric space $X$ of bounded geometry for arbitrarily large $\lambda$
admits a uniformly bounded cover of a finite multiplicity with the Lebesgue 
number $\ge\lambda$. Like in the proof of Lemma 3.1 we can conclude that
every metric space of bounded geometry has the following property:

{(**) \it For arbitrarily large $\lambda$ there is $n=n(\lambda)$, $n$-dimensional
uniform simplicial complex $K_{\lambda}$ of size $\lambda$, supplied with
$l_{\infty}$-metric, and a 1-Lipschitz uniformly cobounded map
$\phi_{\lambda}:X\to K_{\lambda}$}.

A map of a metric space to a simplicial complex $f:X\to K$ is called uniformly
cobounded if there is a constant $C$ such that $\diam \phi^{-1}(\Delta)\le  C$ 
for all simplices $\Delta$. The property (**) can be organized in so called
an anti-\v Cech approximation of $X$ [H-R1]. An {\it anti-\v Cech approximation} of 
$X$ is a direct
system of finite dimensional uniform simplicial complexes with 1-Lipschitz 
bonding maps $\{ K_{\lambda_n},g^m_n\}$ with $\lambda_n\to\infty$ together 
with a sequence of 1-Lipschitz uniformly cobounded maps 
$\phi_{\lambda_n}:X\to K_{\lambda_n}$ such that 
$\phi_{\lambda_n}=g^m_n\phi_{\lambda_m}$ for all
$n>m$. We call an anti-\v Cech approximation {\it Euclidean} if all complexes
$K_{\lambda_n}$ are given the Euclidean metric.

We say that $X$ admits a \v Cech approximation by a certain class $\sC$ of asymptotic
polyhedra if given a proper function $f:X\to\R_+$ there is an asymptotic
polyhedron $N\in\sC$ and an 1-Lipschitz map $\phi:X\to N$ with
$\diam \phi^{-1}(\Delta)\le\min_{x\in\Delta}f(x)$ for every simplex 
$\Delta\subset N\setminus C$ for some compact set $C$. We recall that Lemma 3.1
asserts that if $\as X<\infty$, then $X$ can be approximated by
$n$-dimensional Euclidean asymptotic polyhedra $K$ with the following properties:

1) Every simplex in $K$ is either the standard $\Delta_{2^i}$ of size $2^i$ for 
some $i$ or
it isomorphic to the simplex $\tilde\Delta_{2^i}$ with edges of two types: 
of  length $2^i$ and $2^{i+1}$; 
 
2) If two simplices in $K$ have nonempty intersection
then they either of the same type or they have the types $\Delta_{2^i}$ and
$\tilde\Delta_{2^i}$ or they have the types $\Delta_{2^i}$ and 
$\tilde\Delta_{2^{i-1}}$.

The class of Euclidean asymptotic polyhedra with the above property we denote
by $\sC_0$.
The arguments of Lemma 3.1 and Lemma 3.2, applied in the infinite dimensional
case, gives us the following
\proclaim{Proposition 6.1}
If a metric space of bounded geometry $X$ admits an Euclidean 
anti-\v Cech approximation, then it can be approximated by Euclidean asymptotic 
polyhedra from the class $\sC_0$.
\endproclaim

\proclaim{Proposition 6.2}
If a metric space of bounded geometry $X$ admits an Euclidean 
anti-\v Cech approximation, then it
admits an approximation by Euclidean asymptotic polyhedra $\phi:X\to K$ from the 
class
$\sC_0$ with proper homotopy inverse maps $\gamma:X\to K$.
\endproclaim

DEFINITION. A metric space that admits a  \v Cech) approximation by Euclidean
asymptotic polyhedra
is called {\it asymptotically piecewise Euclidean}.

We note that if $X$ admits an approximation by Euclidean asymptotic polyhedra
then it admits an approximation by asymptotic polyhedra from the class $\sC_0$.

We recall that a {\it uniform embedding} in the coarse sense is
a map $q:X\to Y$ between metric spaces such that
$$\rho_1(d_X(x,y))\le d_Y(q(x),q(y))\le\rho_2(d_X(x,y))$$ for two 
tending to infinity functions $\rho_1,\rho_2:\R_+\to\R_+$
and all $x,y\in X$.
\proclaim{Lemma 6.3}
Every Euclidean asymptotic polyhedron $K\in\sC_0$ admits a uniform embedding 
in the Hilbert space $l_2$.
\endproclaim
\demo{Proof}
The uniform simplicial complex $K_U$ can be naturally embedded into the
infinite simplex $\Delta^{\infty}=\{(t_1,t_2,\dots)\in l_2\mid \Sigma t_i=1,
t_i\ge 0\}$. 
Let $\Delta$ be a simplex in $K$ and let
$l(\Delta)$ denote the minimum of the length of edges in $\Delta$.
We define a map $l:K\to\R_+$ by the rule:
$l(x)=\min\{l(\Delta)\mid x\in\Delta\}$.
Let $\rho:K\to\R_+$ be an 1-Lipschitz function tending to infinity and
$l/2\le\rho\le l$. We define a map $q:K\to l_2$ by the formula 
$q(x)=\rho(x)u(x)$ and show that it is a uniform embedding. 

First we show that $q$ is 3-Lipschitz. 
Since a metric on $K$ is geodesic, it suffices to show that $q$ is
3-Lipschitz on every simplex.
For every two points $x$ and $y$ from the same
simplex $\Delta$ we have
$$\|q(x)-q(y)\|=\|\rho(x)u(x)-\rho(y)u(y)\|
\le \|\rho(x)u(y)-\rho(y)u(y)\|+\|\rho(x)u(x)-\rho(x)u(y)\|\le$$
$$d_K(x,y)\|u(y)\|+\rho(x)\|u(x)-u(y)\|\le d_K(x,y)+\rho(x) d_U(x,y).$$
Since the projection of the simplices $\Delta_l$ and $\tilde\Delta_l$
onto the standard unit simplex
is $(1/l)$-Lipschitz, we have $$d_U(x,y)\le (2/l(\Delta))d_K(x,y)\le 
(2/l(x))d_K(x,y)\le(2/\rho(x))d_K(x,y).$$
Then $\|q(x)-q(y)\|\le 3d_K(x,y)$.

Now assume that $q^{-1}$ is not uniform. It means that there is a number $c>0$
and a sequence of pairs of points $\{(x_k,y_k)\}$ in $K$ such that
$d_K(x_k,y_k)\to\infty$ and $\|q(x_k)-q(y_k)\|\le c$. We may assume that
$\rho(x_k)\le\rho(y_k)$ for all $k$ . Then $c/\rho(x_k)\ge\|u(x_k)-u(y_k)\|$.
Since $\rho(x_k)\to\infty$, for large enough $k$ we have the
inequality $\|u(x_k)-u(y_k)\|\ge(1/2)d_U(x_k,y_k)$. Moreover, we may assume
that $x_k$ and $y_k$ lie in two simplices $\Delta_1$, $\Delta_2$ in $K$ with
$\Delta_1\cap\Delta_2\ne\emptyset$. Then
 $2\max\{l(\Delta_1),l(\Delta_2)\}d_U(x_k,y_k)\ge d_K(x_k,y_k)$.
Therefore, $4\max\{l(x),l(y)\}d_U(x_k,y_k)\ge d_K(x_k,y_k)$. Then
$16\rho(x)d_U(x_k,y_k)\ge d_K(x_k,y_k)$ and hence we get a contradiction:
$32c\ge 16\rho(x)d_U(x_k,y_k)\ge d_K(x_k,y_k)$ for all $k$. 
\qed
\enddemo

Clearly, every discrete subspace of the Hilbert space is asymptotically 
Euclidean.
\proclaim{Theorem 6.4}
For an asymptotically Euclidean geometrically finite
groups $\Gamma$ the Roe index map is a split monomorphism.
\endproclaim
\demo{Proof}
Let $X=E\Gamma$ 
 and let $\phi:X\to K$ and $\gamma: K\to X$ as in Proposition 6.2.
According to Lemma 6.3 the asymptotic polyhedron $K$ is uniformly embeddable 
into $l_2$.
By the theorem of Yu [Yu2] the coarse Baum-Connes conjecture holds for $K$, i.e.
the Roe index map $A_K:K_*^(K)\to K_*(C*K)$ is an isomorphism.
Then the following diagram implies the splitting of $A_X$:
$$
\CD
K_*(X) @>A_X>> K_*(C^*X)\\
@V\phi_*VV @V\phi'VV\\
K_*(K) @>A_K>> K_*(C^*K)\\
@V\gamma_*VV\\
K_*(X)\\
\endCD
$$
\qed
\enddemo
It is known [Ro1] that the monomorphism version of the coarse Baum-Connes
conjecture implies the Gromov-Lawson conjecture. This together with Theorem 6.4
gives the following
\proclaim{Corollary 6.5}
The Gromov-Lawson conjecture holds for manifolds with 
asymptotically piecewise Euclidean universal covers. 
\endproclaim

In view of the main result of the next section an affirmative answer to the 
following question would eliminate an approach to the Novikov Conjecture
via hypersphericity (see a remark on page 35 of [G4]).
\proclaim{Problem 6.6}
Is every hypereuclidean (hyperspherical) manifold asymptotically piecewise Euclidean?
\endproclaim

\head \S 7 Expanders are not asymptotically piecewise Euclidean\endhead

Let $X$ be a finite graph, we denote by $V$ the set of vertices and by $E$
the set of edges in $X$. We will identify the graph $X$ with its set of vertices 
$V$. Every graph is a metric space with respect to the natural metric where every 
edge has the length one. For a subset $A\subset X$ we define the boundary
$\partial A=\{x\in X\mid dist(x,A)=1\}$. Let $|A|$ denote the 
cardinality of $A$.

DEFINITION [Lu]. An {\it expander} with a conductance number $c$ and the degree 
$d$ is an infinite sequence of finite graphs $\{X_n\}$ with the degree $d$ such that 
$|X_n|$ tends to infinity and for every $A\subset X_n$ with $|A|\le |X_n|/2$
there is the inequality $|\partial A|\ge c|A|$.

Let $X$ be a finite graph, we denote by $P$ all nonordered pairs of distinct
points in $X$. For every nonconstant map $f:X\to l_2$ to the Hilbert space
we introduce the number 
$$
D^2_f=
\frac{\frac{1}{|K|}\Sigma_{\{x,y\}\in P}\|f(x)-f(y)\|^2}{\frac{1}{|E|}
\Sigma_{\{x,y\}\in E}\|f(x)-f(y)\|^2}.
$$
If $X$ is a graph with the degree $d$ and with $|X|=n$, then $|P|=n(n-1)/2$
and $|E|=dn/2$.

The following Lemma is well-known. It can be derived from [M, Proposition 3].
It also can be obtain from the equality
$$
\lambda_1(X)=\inf\{\frac{\|df\|^2}{\|f\|^2}\mid \Sigma f(x)=0\}
$$
for the first positive eigenvalue of the Laplacian on $X$ and the Cheeger's 
inequality (see Proposition 4.2.3 in [Lu]). 
\proclaim{Lemma 7.1} Let $\{X_n\}$ be an expander. Then there is a constant
$c_0$ such that $D^2_{f_n}\le c_0$ for all $n$ for all possible maps
 $f_n:X_n\to l_2$  to the Hilbert space.
\endproclaim
\proclaim{Corollary 7.2}
For every sequence of 1-Lipschitz maps $f_n:X_n\to l_2$ there is the inequality
$\frac{1}{|P_n|}\Sigma_{P_n}\|f(x)-f(y)\|^2\le c_0$ for every $n$.
\endproclaim
\demo{Proof}
In the case of 1-Lipschitz map we have 
$${\frac{1}{|E|}\Sigma_{\{x,y\}\in E}\|f(x)-f(y)\|^2}\le 1.$$\qed
\enddemo

We say that a metric space $X$ contains an expander $\{X_n\}$ if there is 
a sequence of isometric embeddings $X_n\to X$.
\proclaim{Theorem 7.3}
An asymptotically piecewise Euclidean metric space $X$ cannot contain an expander.
\endproclaim
\demo{Proof}
Assume the contrary, i.e. assume that $X$ contains an expander of degree $d$.
We enumerate all the graphs in the expander by a subsequence of $\N$ in such a 
way that $|X_n|=n$. 
Let $x_0\in X$ be a base point.
Let $n(t)=\min\{n\mid
X_n\cap\partial B_t(x_0)\ne\emptyset\}$. We note that $n(t)$ tends to infinity
as $t$ approaches infinity. We consider a proper function $f:X\to\R_+$
such that $f(x)<\log_dn(\|x\|)/4$. Then like in the proof of Lemma 3.1
one can take an 1-Lipschitz map $\phi:X\to N$ to an asymptotic polyhedron such that
for every $r>0$ there is a compact set $C_r\subset X$ with
$diam(\phi^{-1}(B_r(y)))\le f(z)$ for every $z\in\phi^{-1}(B_r(y))\setminus C_r$.
Since $X$ is asymptotically Euclidean, we may assume that $N$ is an Euclidean
asymptotic polyhedron. By Lemma 6.3 there is a uniform embedding
$g:N\to l_2$. By the definition, there is a monotone function $\rho:\R_+\to\R_+$
tending to infinity such that $$\rho(d_N(x',y'))\le\|g(x')-g(y')\|\le d_N(x',y')$$ 
for all $x',y'\in N$.

We take $r$ such that $\rho(r)>2\sqrt{c_0}$ and
we take $n$ sufficiently large.
Denote by $P'_n\subset P_n$ the set of pairs $\{x,y\}$ in $X_n\setminus C_r$ with
$d(x,y)\ge\log_d(n/4)$. Since the degree of $X_n$ is $d$, we have 
$$|B_k(v)\cap X_n|\le 1+d+d^2+\dots+d^k\le 2d^k.$$
Therefore at least $n(n-2\frac{n}{4})/2=n^2/4$ such pairs are contained in $X_n$. 
Then for
sufficiently large $n$ we have at least $n^2/8$ such pairs in $X_n\setminus C_r$. 
Show that every such pair $\{x,y\}$ in $X_n\setminus C_r$ satisfies the
inequality $d_N(\phi(x),\phi(y))>r$. Indeed, if we assume that
$d_N(\phi(x),\phi(y))\le r$, then we will obtain that
$$diam\phi^{-1}(B_r(\phi(y)))\ge d(x,y)\ge\log_d(n/4).$$ On the other hand, we have 
that $$diam\phi^{-1}(B_r(\phi(y)))\le f(x)<\log_d(n(\|x\|)/4).$$ Since
$n\ge n(\|x\|)$, we obtain the contradiction 
$$diam\phi^{-1}(B_r(\phi(y)))<\log_d(n/4).$$

Let $f=g\circ\phi$. By the Corollary 7.2 we have
$$
c_0\ge\frac{1}{|P_n|}\Sigma_{P_n}\|f(x)-f(y)\|^2\ge
\frac{1}{|P_n|}\Sigma_{P'_n}\|f(x)-f(y)\|^2\ge\frac{n^2}{8|P_n|}\min_{P'_n}
\|f(x)-f(y)\|^2=$$

$$\frac{n^2}{4n(n-1)}\|f(\tilde x)-f(\tilde y)\|^2\ge\frac{1}{4}\|f(\tilde x)-
f(\tilde y)\|^2=\frac{1}{4}\|g(\phi(\tilde x))-g(\phi(\tilde y))\|^2\ge$$

$$\frac{1}{4}\rho^2(d_N(\phi(\tilde x),\phi(\tilde y)))\ge\frac{1}{4}\rho^2(r)>c_0.$$
 This contradiction completes the proof.
\qed
\enddemo

REMARK. One can show (see [M]) that for every $p\ge 1$ the number
$$
D^p_f=\frac{\frac{1}{|P|}\Sigma_{\{x,y\}\in K}\|f(x)-f(y)\|^p}
{\frac{1}{|E|}\Sigma_{\{x,y\}\in E}\|f(x)-f(y)\|^p}\le c_0
$$
is bounded from above by the same number $c_0$ for all maps $f:X\to l_p$ to
the Banach space $l_p$. This implies in particular that an expander is not uniformly
embeddable in $l_p$ for any $p$.
\proclaim{Problem 7.4}
Given $p>0$, does the Novikov conjecture holds for a group $\Gamma$ which
admits a uniform embedding in $l_p$ ?
\endproclaim
It is known that the answer is `yes' for $p=2$ [Yu2] and for $p=1$ [D-G-L-Y].

\enddocument